\documentclass{article}
\usepackage{amssymb}
\usepackage{graphicx}
\usepackage{amsmath}

\input{tcilatex}
\begin{document}

\begin{center}
\textbf{Nonlocal fractional elliptic and parabol\i c equat\i ons in Besov
spaces and applications\ \ }

\textbf{Veli Shakhmurov }

Antalya Bilim University, \c{C}\i plakl\i\ Mah. Farabi Cad. 23 Dosemealti
07190 Antalya, Turkey, E-mail: veli.sahmurov@gmail.com

Baku Engineering University, 120 Hasan Aliyev Absheron, Baku\ 0101 Azerbaijan

\ \ \ \ \ \ \ \ \ 
\end{center}

\begin{quote}
\ \ \ \ \ \ \ \ \ \ \ \ \ \ \ \ \ \ \ \ \ \ \ \ \ \ \ \ \ \ \ \ \ \ \ \ \ \
\ \ \ \ \ \ \ \ \ \ \ \ \ \ \ 
\end{quote}

\begin{center}
\textbf{ABSTRACT}
\end{center}

\begin{quote}
\ \ \ \ \ \ \ \ \ \ \ \ \ \ \ \ \ 
\end{quote}

\ \ The maximal $B_{p,q}^{s}$-regularity properties of a nonlocal fractional
elliptic equation is studied. Particularly, it is proven that the operator
generated by this nonlocal ell\i p\i t\i c equation  in $B_{p,q}^{s}$ is
sectorial and also is a generator of an analytic semigroup. Moreover,
well-posedeness of\ nonlocal fractional parabolic equation in Besov spaces
are established.

\begin{center}
\bigskip\ \ \textbf{AMS: 47GXX, 35JXX, 47FXX, 47DXX, 43AXX}
\end{center}


\begin{center}
\textbf{1. Introduction, notations and background }
\end{center}

In the last years, fractional elliptic and parabolic equations have found
many applications in physics (see $\left[ \text{2, 5}\right] $, $\left[ 
\text{7-11}\right] $, $\left[ \text{18}\right] $ and the references
therein). The regularity properties of fractional differential equations
(FDEs) have been studied e.g. in $\left[ \text{3, 8, 9, 12-16}\right] $. The
main objective of the present paper is to discuss the $B_{p,q}^{s}\left( 
\mathbb{R}^{n}\right) $-maximal regularity of the nonlocal elliptic FDE with
parameter%
\begin{equation}
\sum\limits_{\left\vert \alpha \right\vert \leq l}a_{\alpha }\ast D^{\alpha
}u+\lambda u=f\left( x\right) \text{, }x\in \mathbb{R}^{n},  \tag{1.1}
\end{equation}%
where $a_{\alpha }$ are complex valued functions, $\lambda $ is a complex
parameter and%
\begin{equation*}
D^{\alpha }=D_{1}^{\alpha _{1}}D_{2}^{\alpha _{2}}...D_{n}^{\alpha _{n}}%
\text{,}\ \alpha _{i}\in \left[ 0,\infty \right) \text{, }\alpha =\left(
\alpha _{1},\alpha _{2},...,\alpha _{n}\right) .
\end{equation*}%
Here $D_{k}^{\alpha _{k}}$are Caputo type fractional partial derivatives of
order $\alpha _{k}\in \left[ {}\right. m-1,m\left. {}\right) $ with respect
to $x_{k}\in \left( a,b\right) $ i.e. 
\begin{equation}
D_{k}^{\alpha _{k}}u=\frac{1}{\Gamma \left( m-\alpha _{k}\right) }%
\dint\limits_{a}^{x_{k}}\left( x_{k}-\tau \right) ^{m-\alpha
_{k}-1}u^{\left( m\right) }\left( \tau \right) d\tau \text{,}  \tag{1.2}
\end{equation}
$\Gamma \left( \alpha _{k}\right) $ is Gamma function for $\alpha _{k}>0$
(see e.g. $\left[ \text{5, 7, 11}\right] $), $m\geq 1$ is a positive
integer, $a_{\alpha }\ast u$ is a convolution of \ $a_{\alpha }$ and $u$.

Let $E$ be a Banach space. Here, $L_{p}\left( \Omega \text{;}E\right) $
denotes the space of $E$-valued strongly measurable complex-valued functions
that are defined on the measurable subset $\Omega \subset \mathbb{R}^{n}$
with the norm given by

\begin{equation*}
\left\Vert f\right\Vert _{L_{p}\left( \Omega \text{;}E\right) }=\left(
\int\limits_{\Omega }\left\Vert f\left( x\right) \right\Vert
_{E}^{p}dx\right) ^{\frac{1}{p}}\text{, }1\leq p<\infty \ .
\end{equation*}

Let $m_{i}$, $s_{i}$ be positive integers, $k_{i}$ be nonnegative integers, $%
m_{i}>s_{i}-k_{i}>0$, $i=1,2,..n$ and $s=\left( s_{1},s_{2},...,s_{n}\right) 
$, $1\leq p\leq \infty $, $1\leq q\leq \infty $, $0<h_{0}<\infty $. The $E$%
-valued Besov space $B_{p,q}^{s}\left( \Omega \text{;}E\right) $ are defined
as 
\begin{equation*}
\begin{array}{c}
B_{p,q}^{s}\left( \Omega \text{;}E\right) =\left\{ f\right. :f\in
L_{p}\left( \Omega \text{;}E\right) ,\left\Vert f\right\Vert
_{B_{p,q}^{s}\left( \Omega \text{;}E\right) }= \\ 
\left\Vert f\right\Vert _{L_{p}\left( \Omega \text{;}E\right)
}+\sum\limits_{i=1}^{n}\left( \int\limits_{0}^{h_{0}}y^{-\left[ \left(
s_{i}-k_{i}\right) q+1\right] }\left\Vert \Delta _{i}^{m_{i}}\left( y,\Omega
\right) D_{i}^{k_{i}}f\right\Vert _{L_{p}\left( \Omega \text{;}E\right)
}^{q}dy\right) ^{\frac{1}{q}}<\infty ,%
\end{array}%
\end{equation*}%
\ \ 

\ 
\begin{equation*}
\left\Vert f\right\Vert _{B_{p,q}^{s}\left( \Omega \text{;}E\right)
}=\sum\limits_{i=1}^{n}\sup\limits_{0<y<y_{0}}\frac{\left\Vert \Delta
_{i}^{m_{i}}\left( y,\Omega \right) D_{i}^{k_{i}}f\right\Vert _{L_{p}\left(
\Omega \text{;}E\right) }}{y^{s_{i}-k_{i}}}\text{ \ for }\theta =\infty .
\end{equation*}

Let $E_{0}$ and $E$ be two Banach spaces and $E_{0}$ be continuously and
densely embedded into $E$. Moreover, let $l\in \mathbb{R}_{+}$. Consider the
Sobolev-Besov space $B_{p,q}^{l,s}\left( \Omega ;E_{0},E\right)
=B_{p,q}^{l,s}\left( \Omega ;E\right) \cap B_{p,q}^{s}\left( \Omega
;E_{0}\right) $ with the norm%
\begin{equation*}
\left\Vert u\right\Vert _{B_{p,q}^{l,s}\left( \Omega ;E_{0},E\right)
}=\left\Vert u\right\Vert _{B_{p,q}^{s}\left( \Omega ;E\right)
}+\dsum\limits_{\left\vert \alpha \right\vert \leq l}\left\Vert D^{\alpha
}u\right\Vert _{B_{p,q}^{s}\left( \Omega :E\right) }<\infty .
\end{equation*}

Here, $\mathbb{C}$ denotes the set of complex numbers. For $E_{0}=E=\mathbb{C%
}$ the spaces $L_{p}\left( \Omega \text{;}E\right) $, $B_{p,q}^{s}\left(
\Omega \text{;}E\right) $, $B_{p,q}^{l,s}\left( \Omega ;E_{0},E\right) $
will be denoted by $L_{p}\left( \Omega \right) $, $B_{p,q}^{s}\left( \Omega
\right) $ and\ $B_{p,q}^{l,s}\left( \Omega \right) $, respectively.

Let $S\left( \mathbb{R}^{n}\right) $ denote Schwartz class, i.e., the space
of all rapidly decreasing smooth functions on $\mathbb{R}^{n}$ equipped with
its usual topology generated by seminorms. A function $\Psi \in C\left( 
\mathbb{R}^{n}\right) $ is called a Fourier multiplier from $%
B_{p,q}^{s}\left( \mathbb{R}^{n}\right) $\ to $B_{p,q}^{s}\left( \mathbb{R}%
^{n}\right) $ if the map%
\begin{equation*}
u\rightarrow \Lambda u=F^{-1}\Psi \left( \xi \right) Fu\text{, }u\in S\left( 
\mathbb{R}^{n}\right)
\end{equation*}%
is well defined and extends to a bounded linear operator 
\begin{equation*}
\Lambda :\ B_{p,q}^{s}\left( \mathbb{R}^{n}\right) \rightarrow \
B_{p,q}^{s}\left( \mathbb{R}^{n}\right) .
\end{equation*}

We prove that problem $\left( 1.1\right) $ has a unique solution $u\in $\ $%
B_{p,q}^{l,s}\left( \mathbb{R}^{n}\right) $ for $f\in $ $B_{p,q}^{s}\left( 
\mathbb{R}^{n}\right) $ and the following uniform coercive estimate holds 
\begin{equation}
\dsum\limits_{\left\vert \alpha \right\vert \leq l}\left\vert \lambda
\right\vert ^{1-\frac{\left\vert \alpha \right\vert }{l}}\left\Vert
a_{\alpha }\ast D^{\alpha }u\right\Vert _{B_{p,q}^{s}\left( \mathbb{R}%
^{n}\right) }\leq C\left\Vert f\right\Vert _{B_{p,q}^{s}\left( \mathbb{R}%
^{n}\right) }.  \tag{1.3}
\end{equation}%
\ The estimate $\left( 1.3\right) $ implies that the operator $O$ that
generated by problem $\left( 1.1\right) $ has a bounded inverse from $%
B_{p,q}^{s}\left( \mathbb{R}^{n}\right) $ into the space $%
B_{p,q}^{l,s}\left( \mathbb{R}^{n}\right) $. Particularly, from the estimate 
$\left( 1.3\right) $ we obtain that $O$ is uniformly sectorial operator in $%
B_{p,q}^{s}\left( \mathbb{R}^{n}\right) $. By using the coercive properties
of elliptic operator, we prove the well posedness of the Cauchy problem for
the nonlocal fractional parabolic differential equation:

\begin{equation}
\frac{\partial u}{\partial t}+\sum\limits_{\left\vert \alpha \right\vert
\leq l}a_{\alpha }\ast D^{\alpha }u=f\left( t,x\right) \text{, }u(0,x)=0, 
\tag{1.4}
\end{equation}%
in the Besov space 
\begin{equation*}
Y^{s}=B_{p_{1},q}^{s}\left( \mathbb{R}_{+};B_{p,q}^{s}\left( \mathbb{R}%
^{n}\right) \right) .
\end{equation*}%
In other words, we show that problem $\left( 1.4\right) $ for each $f\in
Y^{s}$ has a unique solution $u\in B_{\mathbf{p}\text{,}q}^{1,l,s}\left( 
\mathbb{R}_{+}^{n+1}\right) $ with $\mathbf{p=}\left( p,p_{1}\right) $\
satisfying the coercive estimate 
\begin{equation}
\left\Vert \frac{\partial u}{\partial t}\right\Vert
_{Y^{s}}+\sum\limits_{\left\vert \alpha \right\vert \leq l}\left\Vert
a_{\alpha }\ast D_{x}^{\alpha }u\right\Vert _{Y^{s}}+\left\Vert A\ast
u\right\Vert _{Y^{s}}\leq M\left\Vert f\right\Vert _{Y^{s}}.  \tag{1.5}
\end{equation}%
Let 
\begin{equation*}
S_{\varphi }=\left\{ \lambda \text{: \ }\lambda \in \mathbb{C}\text{, }%
\left\vert \arg \lambda \right\vert \leq \varphi \right\} \cup \left\{
0\right\} \text{, }0\leq \varphi <\pi .\ 
\end{equation*}

$L\left( E_{1},E_{2}\right) $ denotes the space of bounded linear operators
from $E_{1}$ to $E_{2}$. For $E_{1}=E_{2}=E$ it denotes by $L\left( E\right) 
$. Let $D\left( A\right) $, $R\left( A\right) $ denote the domain and range
of the linear operator in $E,$ respectively. Let Ker $A$ denote a null space
of $A$. A closed linear operator\ $A$ is said to be $\varphi $-sectorial (or
sectorial for $\varphi =0$) in a Banach\ space $E$ with bound $M>0$ if Ker $%
A=\left\{ 0\right\} $, $D\left( A\right) $ and $R\left( A\right) $ are dense
on $E,$ and $\left\Vert \left( A+\lambda I\right) ^{-1}\right\Vert _{L\left(
E\right) }\leq M\left\vert \lambda \right\vert ^{-1}$ for all $\lambda \in
S_{\varphi },$ $\varphi \in \left[ 0,\right. \left. \pi \right) $, where $I$
is an identity operator in $E.$ Sometimes $A+\lambda I$\ will be written as $%
A+\lambda $ and will be denoted by $A_{\lambda }$. It is known $\left[ \text{%
17, \S 1.15.1}\right] $ that the powers\ $A^{\theta }$, $\theta \in \left(
-\infty ,\infty \right) $ for a sectorial operator $A$ exist.

Here, $S^{\prime }=S^{\prime }\left( \mathbb{R}^{n}\right) $ denotes the
space of linear continuous mappings from $S\left( \mathbb{R}^{n}\right) $
into\ $\mathbb{C}$ and it is called the Schwartz distributions. For any $%
\alpha =\left( \alpha _{1},\alpha _{2},...,\alpha _{n}\right) $, $\alpha
_{i}\in \left[ 0,\infty \right) $, $\xi =\left( \xi _{1},\xi _{2},...,\xi
_{n}\right) \in \mathbb{R}^{n}$ the function $\left( i\xi \right) ^{\alpha }$
will be defined as:

{\large 
\begin{equation*}
\left( i\xi \right) ^{\alpha }=\left\{ 
\begin{array}{c}
\left( i\xi _{1}\right) ^{\alpha _{1}}...\left( i\xi _{n}\right) ^{\alpha
_{n}}\text{, }\xi _{1}\xi _{2}...\xi _{n}\neq 0 \\ 
0\text{, }\xi _{1},\xi _{2}...\xi _{n}=0,%
\end{array}%
\right.
\end{equation*}%
}where 
\begin{equation*}
\left( i\xi _{k}\right) ^{\alpha _{k}}=\exp \left[ \alpha _{k}\left( \ln
\left\vert \xi _{k}\right\vert +i\frac{\pi }{2}\text{ sgn }\xi _{k}\right) %
\right] \text{, }k=1,2,...,n.
\end{equation*}

Sometimes we use one and the same symbol $C$ without distinction in order to
denote positive constants which may differ from each other even in a single
context. When we want to specify the dependence of such a constant on a
parameter, say $\alpha $, we write $C_{\alpha }$.

The embedding theorems in vector valued spaces play a key role in the theory
of DOEs. From $\left[ 15\right] $ we obtain the estimating lower order
derivatives

{\large \ }\textbf{Theorem A}$_{1}$\textbf{.} Suppose $1<p\leq p_{1}<\infty $%
, $l$ is a positive integer and $s$ $\in \left( 0,\infty \right) $ with $%
\varkappa =\frac{1}{l}\left[ \left\vert \alpha \right\vert +n\left( \frac{1}{%
p}-\frac{1}{p_{1}}\right) \right] \leq 1$ for $0\leq \mu \leq 1-\varkappa ,$
then the embedding 
\begin{equation*}
D^{\alpha }B_{p,q}^{s,l}\left( \mathbb{R}^{n}\right) \subset
B_{p_{1},q}^{s}\left( \mathbb{R}^{n}\right) 
\end{equation*}%
is continuous and there exists a constant \ $C_{\mu }$ \ $>0$, depending
only on $\mu $ such that 
\begin{equation*}
\left\Vert D^{\alpha }u\right\Vert _{B_{p_{1},q}^{s}\left( \mathbb{R}%
^{n}\right) }\leq C_{\mu }\left[ h^{\mu }\left\Vert u\right\Vert
_{B_{p,q}^{s,l}\left( \mathbb{R}^{n}\right) }+h^{-\left( 1-\mu \right)
}\left\Vert u\right\Vert _{B_{p,q}^{s}\left( \mathbb{R}^{n}\right) }\right] 
\end{equation*}%
for all $u\in B_{p,q}^{s,l}\left( \mathbb{R}^{n}\right) $ and $0<h\leq
h_{0}<\infty .$

\begin{center}
\textbf{2}.\textbf{\ } \textbf{Nonlocal fractional} \textbf{elliptic} 
\textbf{equation }
\end{center}

\bigskip Consider the problem $\left( 1.1\right) $.

\textbf{Condition 2.1. }Assume $a_{\alpha }\in L_{\infty }\left( \mathbb{R}%
^{n}\right) $ such that 
\begin{equation}
L\left( \xi \right) =\sum\limits_{\left\vert \alpha \right\vert \leq l}\hat{a%
}_{\alpha }(\xi )\left( i\xi \right) ^{\alpha }\in S_{\varphi _{1}}\text{, }%
\left\vert L\left( \xi \right) \right\vert \geq
C\sum\limits_{k=1}^{n}\left\vert \hat{a}_{\alpha \left( l,k\right)
}\right\vert \left\vert \xi _{k}\right\vert ^{l},  \tag{2.1}
\end{equation}%
for

\begin{equation*}
\alpha \left( l,k\right) =\left( 0,0,...,l,0,0,...,0\right) \text{, i.e }%
\alpha _{i}=0\text{, }i\neq k.
\end{equation*}

Consider operator functions

\begin{equation}
\sigma _{1}\left( \xi ,\lambda \right) =\lambda \sigma _{0}\left( \xi
,\lambda \right) \text{, }\sigma _{2}\left( \xi ,\lambda \right)
=\sum\limits_{\left\vert \alpha \right\vert \leq l}\left\vert \lambda
\right\vert ^{1-\frac{\left\vert \alpha \right\vert }{l}}\hat{a}_{\alpha
}(\xi )\left( i\xi \right) ^{\alpha }\sigma _{0}\left( \xi ,\lambda \right) ,
\tag{2.2.}
\end{equation}%
where 
\begin{equation*}
\sigma _{0}\left( \xi ,\lambda \right) =\left[ L\left( \xi \right) +\lambda %
\right] ^{-1}\text{ }.
\end{equation*}%
Let 
\begin{equation*}
X=B_{p,q}^{s}\left( \mathbb{R}^{n}\right) \text{, }Y=B_{p,q}^{l,s}\left( 
\mathbb{R}^{n}\right) .
\end{equation*}%
\ In this section we prove the following:

\textbf{Theorem 2.1.} Assume that the Condition 2.1 is satisfied and $p$, $%
q\in \left[ 1,\infty \right] $.  Suppose that $\gamma \in \left( \left. 1,2%
\right] ,\right. $and $\lambda \in S_{\varphi _{2}}$. Then for $f\in X$, $%
0\leq \varphi _{1}<\pi -\varphi _{2}$ and $\varphi _{1}+\varphi _{2}\leq
\varphi $ there is a unique solution $u$ of the equation $\left( 1.1\right) $
belonging to $Y$ and the following coercive uniform estimate holds 
\begin{equation}
\dsum\limits_{\left\vert \alpha \right\vert \leq l}\left\vert \lambda
\right\vert ^{1-\frac{\left\vert \alpha \right\vert }{l}}\left\Vert a\ast
D^{\alpha }u\right\Vert _{X}+\left\Vert u\right\Vert _{X}\leq C\left\Vert
f\right\Vert _{X}.  \tag{2.3}
\end{equation}

\bigskip For proving of Theorem 2.1 we need the following lemmas:

\bigskip \textbf{Lemma 2.1. }Assume Condition 2.1 holds and $\lambda \in
S_{\varphi _{2}}$with $\varphi _{2}\in \left[ 0,\right. \left. \pi \right) $%
, where $\varphi _{1}+\varphi _{2}<\pi $, then the operator functions $%
\sigma _{i}\left( \xi ,\lambda \right) $ are uniformly bounded, i.e.,%
\begin{equation*}
\left\vert \sigma _{i\varepsilon }\left( \xi ,\lambda \right) \right\vert
\leq C,\text{ }i=0,1,2.
\end{equation*}

\textbf{Proof. }By virtue of $\left[ \text{4, Lemma 2.3}\right] ,$ for $%
L(\xi )\in S_{\varphi _{1}}$, $\lambda \in S_{\varphi _{2}}$ and $\varphi
_{1}+\varphi _{2}<\pi $ there exists a positive constant $C$ such that 
\begin{equation}
\left\vert \lambda +L\left( \xi \right) \right\vert \geq C\left( \left\vert
\lambda \right\vert +\left\vert L\left( \xi \right) \right\vert \right) . 
\tag{2.4}
\end{equation}%
Since\ $L(\xi )\in S_{\varphi _{1}}$ in view of $\ $Condition 2.1 and $%
\left( 2.4\right) $ the function $\sigma _{0}\left( \xi ,\lambda \right) $
is uniformly bounded for all $\xi \in \mathbb{R}^{n}$, $\lambda \in
S_{\varphi _{2}}$, i.e.%
\begin{equation*}
\sigma _{0}\left( \xi ,\lambda \right) \leq (\left\vert \lambda \right\vert
+\left\vert L\left( \xi \right) \right\vert )^{-1}\leq M_{0}.
\end{equation*}%
Moreover,\ we have

\begin{equation*}
\left\vert \sigma _{1}\left( \xi ,\lambda \right) \right\vert \leq
M\left\vert \lambda \right\vert (\left\vert \lambda \right\vert +\left\vert
L\left( \xi \right) \right\vert )^{-1}\leq M_{1}.
\end{equation*}

Next, let us consider $\sigma _{2}.$ It is clear to see that 
\begin{equation}
\left\vert \sigma _{2}\left( \xi ,\lambda \right) \right\vert \leq
C\sum\limits_{\left\vert \alpha \right\vert \leq l}\left\vert \lambda
\right\vert \dprod\limits_{k=1}^{n}\left[ \left\vert \xi \right\vert
\left\vert \lambda \right\vert ^{-\frac{1}{l}}\right] ^{\alpha
_{k}}\left\vert \sigma _{0}\left( \xi ,\lambda \right) \right\vert . 
\tag{2.5}
\end{equation}%
By setting $y_{k}=\left( \left\vert \lambda \right\vert ^{-\frac{1}{l}%
}\left\vert \xi _{k}\right\vert \right) ^{\alpha _{k}}$ in the following
well known inequality 
\begin{equation}
y_{1}^{\alpha _{1}}y_{2}^{\alpha _{2}}...y_{n}^{\alpha _{n}}\leq C\left(
1+\sum\limits_{k=1}^{n}y_{k}^{l}\right) \text{, }y_{k}\geq 0\text{, }%
\left\vert \alpha \right\vert \leq l  \tag{2.6}
\end{equation}%
we get 
\begin{equation*}
\left\Vert \sigma _{2}\left( \xi ,\lambda \right) \right\Vert _{B\left(
E\right) }\leq C\sum\limits_{\left\vert \alpha \right\vert \leq l}\left\vert
\lambda \right\vert \left[ 1+\sum\limits_{k=1}^{n}\left\vert \xi
_{k}\right\vert ^{l}\left\vert \lambda \right\vert ^{-1}\right] \left\vert
\lambda +L\left( \xi \right) \right\vert ^{-1}.
\end{equation*}%
Taking into account the Condition 2.1 and $\left( 2.5\right) -\left(
2.6\right) $ we obtain 
\begin{equation}
\left\vert \sigma _{2}\left( \xi ,\lambda \right) \right\vert \leq C\left(
\left\vert \lambda \right\vert +\sum\limits_{k=1}^{n}\left\vert \xi
_{k}\right\vert ^{l}\right) \left( \left\vert \lambda \right\vert
+\left\vert L\left( \xi \right) \right\vert \right) ^{-1}\leq C.  \notag
\end{equation}

\textbf{Lemma 2.2. }Assume Condition 2.1 holds. Suppose $\hat{a}_{\alpha
}\in C^{\left( n\right) }\left( \mathbb{R}^{n}\right) $ and 
\begin{equation}
\text{ }\left\vert \xi \right\vert ^{\left\vert \beta \right\vert
}\left\vert D^{\beta }\hat{a}_{\alpha }(\xi )\right\vert \leq C_{1}\text{, }%
\beta _{k}\in \left\{ 0,1\right\} ,\text{ }\xi \in \mathbb{R}^{n}\backslash
\left\{ 0\right\} \text{, }0\leq \left\vert \beta \right\vert \leq n, 
\tag{2.7}
\end{equation}

Then, operators $\left\vert \xi \right\vert ^{\left\vert \beta \right\vert
}D_{\xi }^{\beta }\sigma _{i}\left( \xi ,\lambda \right) $, $i=0,1,2$ are
uniformly bounded$\mathbf{.}$

\bigskip \textbf{Proof. }Consider the term $\left\vert \xi \right\vert
^{\left\vert \beta \right\vert }D_{\xi }^{\beta }\sigma _{0}\left( \xi
,\lambda \right) $. By using the Conditin 2.1and the estimates $\left(
2.4\right) -\left( 2.6\right) $, we get 
\begin{equation*}
\left\vert \xi _{k}\right\vert \left\vert D_{\xi _{k}}\sigma _{0}\left( \xi
,\lambda \right) \right\vert \leq
\end{equation*}%
\begin{equation*}
\left[ \left\vert \xi _{k}\right\vert \left\vert \frac{\partial }{\partial
\xi _{k}}\hat{a}_{\alpha }\left( \xi \right) \right\vert +\alpha
_{k}\left\vert \hat{a}_{\alpha }\left( \xi \right) \right\vert \right]
\left\vert \dprod\limits_{k=1}^{n}\left( i\xi _{k}\right) ^{\alpha
_{k}}\right\vert \left\vert \left[ L\left( \xi \right) +\lambda \right]
^{-2}\right\vert <\infty .
\end{equation*}

\bigskip It easy to see that operators $\left\vert \xi \right\vert ^{\beta
}D^{\left\vert \beta \right\vert }\sigma _{0}\left( \xi ,\lambda \right) $
contain the similar terms as in $\left\vert \xi _{k}\right\vert \left\vert
D_{\xi _{k}}\sigma _{0}\left( \xi ,\lambda \right) \right\vert $ for all $%
\beta _{k}\in \left\{ 0,1\right\} $. Hence we get 
\begin{equation*}
\left\vert \xi \right\vert ^{\left\vert \beta \right\vert }\left\vert D_{\xi
}^{\beta }\sigma _{0}\left( \xi ,\lambda \right) \right\vert <\infty .
\end{equation*}

In a similar way, by using the Conditin 2.1and the estimates $\left(
2.4\right) -\left( 2.7\right) $ we obtain 
\begin{equation}
\left\vert \xi \right\vert ^{\left\vert \beta \right\vert }\left\vert D_{\xi
}^{\beta }\sigma _{i}\left( \xi ,\lambda \right) \right\vert <\infty \text{, 
}i=1,2.  \tag{2.8}
\end{equation}

\ \textbf{Proof}.\textbf{\ of Theorem 2.1}. By applying the Fourier
transform to equation $\left( 1.1\right) $\ we get 
\begin{equation}
\hat{u}\left( \xi \right) =\sigma _{0}\left( \xi ,\lambda \right) \hat{f}%
\left( \xi \right) \text{, }\sigma _{0}\left( \xi ,\lambda \right) =\left[
L_{\varepsilon }\left( \xi \right) +\lambda \right] ^{-1}.\text{ }  \tag{2.9}
\end{equation}%
Hence, the solution of $\left( 1.1\right) $ can be represented as $u\left(
x\right) =F^{-1}D_{\varepsilon }\left( \xi ,\lambda \right) \hat{f}$ \ and
by Lemma 2.1 there are positive constants $C_{1}$ and $C_{2}$ such that

\begin{equation*}
C_{1}\left\vert \lambda \right\vert \left\Vert u\right\Vert _{X}\leq
\left\Vert F^{-1}\left[ \lambda \sigma _{0}\left( \xi ,\lambda \right) \hat{f%
}\right] \right\Vert _{X}\leq C_{2}\left\vert \lambda \right\vert \left\Vert
u\right\Vert _{X}\text{, }
\end{equation*}%
\begin{equation}
C_{1}\sum\limits_{\left\vert \alpha \right\vert \leq l}\left\vert \lambda
\right\vert ^{1-\frac{\left\vert \alpha \right\vert }{l}}\left\Vert
a_{\alpha }\ast D^{\alpha }u\right\Vert _{X}\leq \left\Vert F^{-1}\left[
\sigma _{2}\left( \xi ,\lambda \right) \hat{f}\right] \right\Vert _{X}\leq 
\tag{2.10}
\end{equation}%
\begin{equation*}
C_{2}\sum\limits_{\left\vert \alpha \right\vert \leq l}\left\vert \lambda
\right\vert ^{1-\frac{\left\vert \alpha \right\vert }{l}}\left\Vert
a_{\alpha }\ast D^{\alpha }u\right\Vert _{X},
\end{equation*}%
Therefore, it is sufficient to show that the operators$\ \sigma _{i}\left(
\xi ,\lambda \right) $ are multipliers in $X.$ But, by Lemma 2.2 and by
virtue of Fourier multiplier theorem in Besov spaces $B_{p,q}^{s}\left( 
\mathbb{R}^{n}\right) $ (see e.g $\left[ \text{6, Corollary 4.11}\right] $)
we get that $\sigma _{i}\left( \xi ,\lambda \right) $ are multipliers in $X$%
. So, we obtai the assertion.

\bigskip \textbf{Result 2.1.} Theorem 2.1 implies that the operator $O$ is
separable in $X$, i.e. for all $f\in X$ \ there is a unique solution $u\in Y$%
\ of the problem $\left( 1.1\right) $, all terms of equation $\left(
1.1\right) $ are also from $X$ and there are positive constants $C_{1}$ and $%
C_{2}$ so that 
\begin{equation*}
C_{1}\left\Vert Ou\right\Vert _{X}\leq \dsum\limits_{\left\vert \alpha
\right\vert \leq l}\left\Vert a_{\alpha }\ast D^{\alpha }u\right\Vert
_{X}+\left\Vert u\right\Vert _{X}\leq C_{2}\left\Vert Ou\right\Vert _{X}.
\end{equation*}

Indeed, if we put $\lambda =1$ in $\left( 2.3\right) $, by Theorem 2.1 we
get the second inequality. So it is remain to prove the first estimate.\ The
first inequality is equivalent to the following estimate%
\begin{equation*}
\sum\limits_{\left\vert \alpha \right\vert \leq l}\left\Vert F^{-1}\hat{a}%
_{\alpha }\left( i\xi \right) ^{\alpha }\hat{u}\right\Vert _{X}\leq
\sum\limits_{\left\vert \alpha \right\vert \leq l}\left\Vert F^{-1}\hat{a}%
_{\alpha }\left( i\xi \right) ^{\alpha }\sigma _{0}\left( \xi ,\lambda
\right) \hat{f}\left( \xi \right) \right\Vert _{X}.
\end{equation*}

So, it suffices to show that the operator functions%
\begin{equation*}
\sigma _{0}\left( \xi ,\lambda \right) \text{, }\sum\limits_{\left\vert
\alpha \right\vert \leq l}\hat{a}_{\alpha }\left( i\xi \right) ^{\alpha
}\sigma _{0}\left( \xi ,\lambda \right) \text{ }
\end{equation*}%
are uniform Fourier multipliers in $X$. This fact is proved in a similar way
as in the proof of Theorem 2.1.

From Theorem 2.1, we have:

\textbf{Result 2.2. }Assume all conditions of Theorem 2.1 are satisfied.
Then, for all $\lambda \in S_{\varphi }$ the resolvent of operator $O$
exists and the following sharp coercive uniform estimate holds 
\begin{equation}
\dsum\limits_{\left\vert \alpha \right\vert \leq l}\left\vert \lambda
\right\vert ^{1-\frac{\left\vert \alpha \right\vert }{l}}\left\Vert a\ast
D^{\alpha }\left( O+\lambda \right) ^{-1}\right\Vert _{L\left( X\right)
}+\left\Vert \left( O+\lambda \right) ^{-1}\right\Vert _{L\left( X\right)
}\leq C.  \tag{2.11}
\end{equation}

Indeed, we infer from Theorem 2.1 that the operator $O+\lambda $ has a
bounded inverse from $X$ to $Y.$ So, the solution\ $u$\ of the equation $%
\left( 1.1\right) $ can be expressed as $u\left( x\right) =\left( O+\lambda
\right) ^{-1}f$ \ for all $f\in X.$ Then estimate $\left( 2.4\right) $
implies the estimate $\left( 2.11\right) .$

\textbf{Theorem 2.2. } Assume that the Condition 2.1 is satisfied,  $p$, $%
q\in \left[ 1,\infty \right] $ and $\lambda \in S_{\varphi _{2}}$. Then for $%
f\in X$, $0\leq \varphi _{1}<\pi -\varphi _{2}$ and $\varphi _{1}+\varphi
_{2}\leq \varphi $ there is a unique solution $u$ of $\left( 1.1\right) $
belonging to $Y$ and the following coercive uniform estimate holds 
\begin{equation}
\dsum\limits_{\left\vert \alpha \right\vert \leq l}\left\vert \lambda
\right\vert ^{1-\frac{\left\vert \alpha \right\vert }{l}}\left\Vert
D^{\alpha }u\right\Vert _{X}+\left\Vert u\right\Vert _{X}\leq C\left\Vert
f\right\Vert _{X}.  \tag{2.12}
\end{equation}

\textbf{Proof}. The estimate $\left( 2.12\right) $ is derived by reasoning
as in Theorem 2.2.

From Theorem 2.2, we have the following results:

\textbf{Result 2.3.} There are positive constants $C_{1}$ and $C_{2}$ so
that 
\begin{equation}
C_{1}\left\Vert Ou\right\Vert _{X}\leq \dsum\limits_{\left\vert \alpha
\right\vert \leq l}\left\Vert D^{\alpha }u\right\Vert _{X}+\left\Vert
Au\right\Vert _{X}\leq C_{2}\left\Vert Ou\right\Vert _{X}.  \tag{2.13}
\end{equation}

From theorem 2.2. we obtain

\textbf{Result 2.4. }Assume all conditions of Theorem 2.2 hold. Then, for
all $\lambda \in S_{\varphi }$ the resolvent of operator $O$ exists and the
following sharp uniform estimate holds 
\begin{equation}
\dsum\limits_{\left\vert \alpha \right\vert \leq l}\left\vert \lambda
\right\vert ^{1-\frac{\left\vert \alpha \right\vert }{l}}\left\Vert
D^{\alpha }\left( O+\lambda \right) ^{-1}\right\Vert _{L\left( X\right)
}+\left\Vert \left( O+\lambda \right) ^{-1}\right\Vert _{L\left( X\right)
}\leq C.  \tag{2.15}
\end{equation}

\textbf{Result 2.5. }Theorem 2.2 particularly implies that the operator $O$
is sectorial in $X.$ Then the operators $O^{s}$ are generators of analytic
semigroups in $X$ for $s\leq \frac{1}{2}$ (see e.g. $\left[ \text{17, \S %
1.14.5}\right] $)$.$

\begin{center}
\bigskip \textbf{3. The Cauchy problem for fractional parabolic equation }
\end{center}

In this section, we shall consider the following Cauchy problem for the
parabolic FDE 
\begin{equation}
\frac{\partial u}{\partial t}+\sum\limits_{\left\vert \alpha \right\vert
\leq l}a_{\alpha }\ast D^{\alpha }u=f\left( t,x\right) ,\text{ }u(0,x)=0%
\text{, }t\in \mathbb{R}_{+}\text{, }x\in \mathbb{R}^{n},  \tag{3.1}
\end{equation}%
where $a$ is a complex number, $D_{x}^{\alpha }$ is the fractional
derivative in $x$ defined by $\left( 1.2\right) $.

By applying Theorem 2.1 we establish the maximal regularity of the problem $%
\left( 3.1\right) $ in mixed Besov spaces. Let $O$ denote the operator
generated by problem $\left( 1.1\right) $ for $\lambda =0$. Let 
\begin{equation*}
X_{0}=B_{p,q}^{s}\left( \mathbb{R}^{n}\right) \text{, }Y^{s}=B_{p_{1},q}^{s}%
\left( \mathbb{R}_{+};X_{0}\right) \text{, }Y_{0}=B_{\mathbf{p},q}^{s}\left( 
\mathbb{R}_{+}^{n+1}\right) \text{, }
\end{equation*}%
\begin{equation*}
Y^{1,l,s}=B_{\mathbf{p},q}^{1,l,s}\left( \mathbb{R}_{+}^{n+1}\right) \text{, 
}\mathbf{p=}\left( p_{1},p\right) .
\end{equation*}

Let $Y^{1,l,s}$ denotes the space of all functions $u\in X$ possessing the
generalized derivative $D_{t}u=\frac{\partial u}{\partial t}\in Y_{0}$\ and
fractional derivatives $D_{x}^{\alpha }u\in Y_{0}$ for $\left\vert \alpha
\right\vert \leq l$ with the norm 
\begin{equation*}
\ \left\Vert u\right\Vert _{Y^{1,l,s}}=\left\Vert u\right\Vert
_{Y_{0}}+\left\Vert \partial _{t}u\right\Vert
_{Y_{0}}+\sum\limits_{\left\vert \alpha \right\vert \leq l}\left\Vert
D_{x}^{\alpha }u\right\Vert _{Y_{0}},
\end{equation*}%
where $u=u\left( t,x\right) .$

Now, we are ready to state the main result of this section.

\textbf{Theorem 3.1.}\ Assume the conditions of Theorem 2.1 hold for $%
\varphi \in \left( \frac{\pi }{2},\pi \right) $ and $p_{1}$, $p$, $q\in %
\left[ 1,\infty \right] $. Then for $f\in Y$ problem $\left( 3.1\right) $
has a unique solution $u\in Y^{1,l,s}$ satisfying the following unform
coercive estimate 
\begin{equation*}
\left\Vert \frac{\partial u}{\partial t}\right\Vert
_{Y^{s}}+\sum\limits_{\left\vert \alpha \right\vert \leq l}\left\Vert
a_{\alpha }\ast D^{\alpha }u\right\Vert _{Y^{s}}+\left\Vert u\right\Vert
_{Y^{s}}\leq C\left\Vert f\right\Vert _{Y^{s}}.
\end{equation*}

\textbf{Proof.} By definition of $X_{0}$ and by definition of the mixed
space $B_{\mathbf{p,}q}^{s}\left( \mathbb{R}_{+}^{n+1}\right) $ for $\mathbf{%
p=}\left( p\text{, }p_{1}\right) $, we have

\begin{equation}
\left\Vert u\right\Vert _{Y^{s}}=\left\Vert u\right\Vert
_{B_{p_{1},q}^{s}\left( \mathbb{R}_{+};X_{0}\right) }=  \tag{3.2}
\end{equation}%
\begin{equation*}
\left\Vert u\right\Vert _{L_{p_{1}}\left( \mathbb{R}_{+};X_{0}\right)
}+\left( \int\limits_{0}^{y_{0}}y^{-\left[ \left( s-k\right) q+1\right]
}\left\Vert \Delta ^{m_{i}}\left( y,\Omega \right) D^{k}f\right\Vert
_{L_{p_{1}}\left( \mathbb{R}_{+};X_{0}\right) }^{q}dh\right) ^{\frac{1}{q}%
}\geq
\end{equation*}%
\begin{equation*}
\left\Vert u\right\Vert _{B_{\mathbf{p},q}^{s}\left( \mathbb{R}%
_{+}^{n+1}\right) }.
\end{equation*}

Hence, the problem $\left( 3.1\right) $\ can be expressed as the following
Cauchy problem for the parabolic equation 
\begin{equation}
\frac{du}{dt}+Ou\left( t\right) =f\left( t\right) ,\text{ }u\left( 0\right)
=0\text{, }t\in \mathbb{R}_{+}.  \tag{3.3}
\end{equation}

Then, by virtue of $\left[ \text{1, Proposition 8.10}\right] $, we obtain
that for $f\in Y^{s}$ problem $\left( 3.3\right) $ has a unique solution $%
u\in B_{p_{1}\text{,}q}^{1,s}\left( \mathbb{R}_{+};D\left( O\right)
,Y^{s}\right) $ satisfying the following estimate 
\begin{equation*}
\left\Vert \frac{du}{dt}\right\Vert _{B_{p_{1},q}^{s}\left( \mathbb{R}%
_{+};X_{0}\right) }+\left\Vert Ou\right\Vert _{B_{p_{1},q}^{s}\left( \mathbb{%
R}_{+};X_{0}\right) }\leq C\left\Vert f\right\Vert _{B_{p_{1},q}^{s}\left( 
\mathbb{R}_{+};X_{0}\right) }.
\end{equation*}

From the Theorem 2.1, relation $\left( 3.2\right) $ and from the above
estimate we get the assertion$.$

\textbf{References}

\ \ \ \ \ \ \ \ \ \ \ \ \ \ \ \ \ \ \ \ \ \ \ \ \ \ \ \ \ \ \ \ \ \ \ \ \ \
\ \ \ \ \ \ \ \ \ \ \ \ \ \ \ \ \ \ \ \ \ \ \ \ \ \ \ \ \ \ \ \ \ \ \ \ \ \
\ \ \ \ \ \ \ \ \ 

\begin{enumerate}
\item H. Amann, Linear and quasi-linear equations,1, Birkhauser, Basel 1995.

\item D. Baleanu, A. Mousalou, Sh. Rezapour, A new method for investigating
some fractional integro-differential equations involving the Caputo-Fabrizio
derivative. Adv. Differ. Equ. 2017 (51) (2017).

\item P. Clement, G. Gripenberg, S-O. Londen, Regularity properties of
solutions of fractional evolution equations, Evolution equations and their
applications in physical and life sciences, (2019), 235-246.

\item C. Dore, S.Yakubov, Semigroup estimates and non coercive boundary
value problems, Semigroup Forum 60 (2000), 93-121.

\item A. A. Kilbas, H. M. Srivastava, J. J. Trujillo, Theory and
Applications of Fractional Differential Equations. North-Holland Mathematics
Studies, vol. 204. Elsevier, Amsterdam (2006).

\item M. Girardi, L. Weis., Operator-valued Fourier multiplier theorems on
Besov spaces, Math. Nachr., 251(2003), 34--51.

\item K. S. Miller, B. Ross, An Introduction to the Fractional Calculus and
Fractional Differential Equations, JohnWiley \& Sons, New York, NY, USA,
1993.

\item V. Lakshmikantham, A. S.\ Vatsala, Basic theory of fractional
differential equations, Nonlinear Analysis: Theory, Methods \& Applications,
69(8)(2008), 2677-2682.

\item V. Lakshmikantham, J. D.\ Vasundhara, Theory of fractional
differential Equations in a Banach space, European journal pure and applied
mathematics, (1)1 (2008), 38-45.

\item Y-N. Li, H-R. Sun, Integrated fractional resolvent operator function
and fractional abstract Cauchy problem, Abstract and Applied Analysis,
(2014), Article ID 430418.

\item I. Podlubny, Fractional Differential Equations San Diego Academic
Press, 1999.

\item V. B. Shakhmurov, R. V. Shahmurov, Maximal B-regular
integro-differential equations, Chin. Ann. Math. Ser. B, 30B(1), (2008),
39-50

\item V. B. Shakhmurov, R.V. Shahmurov, Sectorial operators with convolution
term, Math. Inequal. Appl., 13 (2) (2010), 387-404.

\item V. B. Shakhmurov, H. K. Musaev, Separability properties of
convolution-differential operator equations in weighted $L_{p}$ spaces,
Appl. and Compt. Math. 14(2) (2015), 221-233.

\item V. B. Shakhmurov, Embedding operators in vector-valued weighted Besov
spaces and applications, Journal of Function spaces and Applications, v.
2012(1), (2012). DOI :http://dx.doi.org/10.1155/2012.

\item V. B. Shakhmurov, Maximal regular abstract elliptic equations and
applications, Siberian Math. Journal, (51) 5, (2010), 935-948.

\item H. Triebel, Interpolation theory, Function spaces, Differential
operators, North-Holland, Amsterdam, 1978.

\item A. Shi A, Yu. Bai, Existence and uniqueness of solution to two-point
boundary value for two-sided fractional differential equations, Applied
Mathematics, 4 (6), (2013), 914-918.
\end{enumerate}

\end{document}